\documentclass[review,hidelinks,onefignum,onetabnum]{siamart250211}
\def\tr{\hbox{\rm tr}}

\usepackage{lipsum}
\usepackage{amsfonts}
\usepackage{graphicx}
\usepackage{epstopdf}
\usepackage{algorithmic}
\usepackage{cite}
\usepackage{lineno}
\nolinenumbers
\ifpdf
  \DeclareGraphicsExtensions{.eps,.pdf,.png,.jpg}
\else
  \DeclareGraphicsExtensions{.eps}
\fi

\newsiamremark{remark}{Remark}
\newsiamremark{example}{Example}
\newsiamremark{hypothesis}{Hypothesis}
\crefname{hypothesis}{Hypothesis}{Hypotheses}
\newsiamthm{claim}{Claim}
\newsiamremark{fact}{Fact}
\crefname{fact}{Fact}{Facts}

\headers{DYNAMIC PROGRAMMING PRINCIPLE }{HAN et al.}

\title{Dynamic programming principle for controlled  systems driven by  fractional Brownian motions and related backward stochastic Riccati equations\thanks{Submitted to the editors DATE.
\funding{This work was supported by National Science Foundation of China (Grant numbers 12471417), National Key R$\&$D Program of China (Grant number 2023YFA1009200) and   an NSERC discovery fund RGPIN 2024-05941 and a centennial fund of University of Alberta.}}}

\author{Yuecai Han\thanks{School of Mathematics, Jilin University, Changchun, 130012, China
  (\email{hanyc@jlu.edu.cn}, \email{yuhangl22@mails.jlu.edu.cn}).}
\and Yaozhong Hu\thanks{Department of Mathematical and Statistical Sciences, University of Alberta at Edmonton,
Edmonton, T6G 2G1, Canada (\email{yaozhong@ualberta.ca}).}
\and Yuhang Li\footnotemark[2]
}

\usepackage{amsopn}

\usepackage{cite}
\ifpdf
\hypersetup{
  pdftitle={S$\Delta$Es and BS$\Delta$Es},
  pdfauthor={D. Doe, P. T. Frank, and J. E. Smith}
}
\fi

\begin{document}

\maketitle

% REQUIRED
\begin{abstract}
This work aims to study the stochastic optimal control problems of controlled stochastic differential systems driven by both Brownian motions and fractional Brownian motions via dynamic programming principle. 
The stochastic Hamilton-Jacobi-Bellman equation is obtained. 
Different from the classical case, the obtained stochastic Hamilton-Jacobi-Bellman equation is a backward stochastic partial differential equation driven by both fractional Brownian motions and standard Brownian motions, whose solution is a pair of random fields. For linear quadratic optimal control problem this backward stochastic partial differential equation is explicitly solved by relating it to  a backward stochastic Riccati equation. 
\end{abstract}

% REQUIRED

% REQUIRED

\section{Introduction}
In this paper, we   solve the following stochastic control problem via dynamic programming principle (DPP, in short) method: minimize the following cost functional over $u$ in certain admissible control domain  $\mathcal{U}[t,T]$:
\begin{align}\label{0}
J\left(u;t,x\right)=\mathbb{E}^{\mathcal{F}_t}\left[\int_t^Tf(s,X_s^{t,x,u},u_s)ds+g(X_T^{t,x,u})\right],
\end{align}
subject to
\begin{align}\label{1}
\left\{\begin{array}{l}
 d X_s^{t,x,u}=b\left(s,X_s^{t,x,u},u_s\right)dt+\sigma\left(s,X_s^{t,x,u}\right)\circ dW^H_s\\
\qquad\qquad\qquad +\gamma\left(s,X_s^{t,x,u},u_s\right)dW_s,\quad s\in[t,T],\\
 X_t^{t,x,u}=x.
\end{array}\right.
\end{align}
Here $W$ is a standard Brownian motion and $W^H$ is a fractional Brownian motion with the underlying Brownian motion
$W$. We shall make  more precise about    the coefficients and the problem in the next section.

It is well known that Hamilton-Jacobi-Bellman  dynamic programming principle and Pontryagin's maximum principle are the two most   important methods to solve   stochastic optimal control problems
when the system is driven by Brownian motion (e.g. $\sigma(s,x)=0$ in \eqref{1}). However, for controlled system such as \eqref{1} driven by both Brownian motion and fractional Brownian  motions, researchers have only  studied the  maximum principle rather than DPP since the involvement of the non-Markovian term $W^H$, due to the potential lack of the Markovian property of the value function. This work is to fill this gap to  obtain  a DPP in solving the optimal control problem $\min J(u; t, x) $ subject to the system \eqref{1}.
%driven by both fractional Brownian motions and standard Brownian motions.
%through DPP in a novel way.

DPP was introduced by Bellman in the 1950s for deterministic optimal control problem. Combining with probability theory and It$\hat{\rm o}$'s formula, DDP for stochastic optimal control problem and related Hamiltonian-Jacobi-Bellman (HJB, in short) equation have been extensively 
investigated. There are many research papers and books on  this topic,  among them we   refer to \cite{yong1999stochastic,fleming2012deterministic} and the references therein. 

To describe the long-range phenomenon in real-world, stochastic differential equations (SDEs, in short) driven by fractional Brownian motion and related control problems have recently been extensively studied. Some related works about SDEs driven by fractional Brownian motion are referred  to \cite{lyons2002system,coutin2002stochastic,rascanu2002differential,hu2007differential,hu2009rough,friz2010multidimensional,hu2013multiple}.  For the control problem, Hu and Zhou\cite{hu2005stochastic} and Duncan and Pasik-Duncan \cite{duncan2013linear} obtained the Markovian control for linear quadratic (LQ, in short) case. Hu and Peng \cite{hu2009backward} investigated backward stochastic differential equations (BSDEs in short) driven by fractional Brownian motions.

For nonlinear control problems driven by fractional Brownian motion,  the  maximum principle has  been largely studied. Han et al.  \cite{han2013maximum} studied the maximum principle for general control system driven by fractional Brownian motion. The adjoint equation   is a  linear BSDEs driven by both fractional Brownian motion and  the underlying standard Brownian motion. Buckdahn and Jing \cite{buckdahn2017mean} investigated the control problem with mean-field. Sun \cite{sun2021stochastic} studied  the control problem of nonconvex   control.   Zheng and Hu \cite{zheng2024global} studied partially observed  stochastic control problem. Han and Li \cite{han2025maximum} studied the control problems driven by fractional noises in discrete-time version. But to the best of our knowledge,  until now there is no  any work concerning  DPP for stochastic optimal control problems driven by fractional Brownian motions.

The objective of this work is to fill out this important gap    to  establish the DPP for optimal control problem (\ref{0})-(\ref{1}). a
Due to the lack of Markovian property, we cannot expect the  HJB equation to be a classical one as in Markovian case. On the contrary, it will be a  (backward) nonlinear  stochastic partial differential equations (BSPDEs, in short) driven by both fractional Brownian motions and standard Brownian motions. As a by product of our general result,  the Feynman-Kac formula for SDEs driven by fractional Brownian motions is also obtained. The relationship between  the  maximum principle and DPP 
%for control problems driven by fractional Brownian motion 
will be established along the way.

As an application  the following LQ control problem is %explicitly 
investigated: minimize the quadratic cost functional 
\begin{align}
J(u;t,x)=\mathbb{E}^{\mathcal{F}_t}\left[\int_t^T\left(S_s|X_s^{t,x,u}|^2+R_su_s^2 \right)ds+G|X_T^{t,x,u}|^2\right],\label{e.1.3}
\end{align}
subject to the following controlled linear system
\begin{align}
\left\{\begin{array}{l}
 d X_s^{t,x,u}=\left(A_sX_s^{t,x,u}+B_su_s\right) dt+(C_sX_s^{t,x,u}+D_su_s)dW_s\\
 \qquad\qquad\qquad+F_sX_s^{t,x,u}\circ dW^H_s,\quad s\in[t,T],\\
 X_t^{t,x,u}=x.
\end{array}\right.\label{e.1.4}
\end{align}
By using the DPP we obtain  that 
the optimal control $u^*$ is given as follows. 
\begin{align}
u^*_s
=&-\big(R_s+D_s^2K_s\big)^{-1}\left((B_s+C_sD_s)K_s+D_s\Theta_s\right)X_s^{t,x,u^*},\label{e.1.5}
\end{align}
where $(K, \Theta)$  is the solution pair to  the following backward stochastic Riccati equation (BSRE, in short):
\begin{equation}\label{2}
    \left\{\begin{array}{l}
 -dK_t=[-(R_t+D_t^2K_t)^{-1}(B_tK_t+D_t\Theta_t+C_tD_tK_t)^2\\
 \quad\qquad\quad\qquad +(2A_t+C_t^2)K_t+2C_t\Theta_t+S_t]dt\\
 \qquad\qquad\qquad +2F_tK_t\circ dW_t^H-\Theta_tdW_t,\quad t\in[0,T],\\
K_T=G.
\end{array}\right.
\end{equation}
Due to the non-Markovian of the LQ  control system, the feedback function
\begin{align*}
-\big(R_s+D_s^2K_s\big)^{-1}\left((B_s+C_sD_s)K_s+D_s\Theta_s\right),
\end{align*}
naturally arisen  from DPP, is an  adapted stochastic process described by BSRE (\ref{2}).  In fact, this type of  BSRE was first introduced to solve control problems with random coefficients in  Tang \cite{tang2003general,tang2015dynamic}.

Let us mention that the optimal Markov control for  one dimensional LQ problem (the space dimension is $1$ in \eqref{e.1.4})  was  solved in \cite{hu2005stochastic} by using the completing square technique.  
However, \eqref{e.1.5}-\eqref{2} tells us that the optimal control is no longer Markovian since $(K_t, \Theta_t)$ depends on $(W_s, W^H_s, 0\le s\le t)$. In fact, $W_s$ (and hence $W^H_s$) can also be represented by $(X_u, 0\le u\le s)$. Thus, $(K_t, \Theta_t)$ depends on $(X_s,  0\le s\le t)$. 

%introduced BSRE driven by standard Brownian motion and prove the existence and uniqueness result. 

The  remaining part   of this paper is organized as follows.  
We obtain  the DPP for optimal control problems driven by fractional Brownian motions and obtain the HJB equations in Section 3 after briefly recall some definitions and results about fractional Brownian motion in Section 2.  We also  establish the   Feynman-Kac formula and explore the relationship between stochastic maximum principle with DPP in Section 3.  In Section 4, the obtained DPP is applied to  LQ control problem  and a BSRE is obtained. In Section 5, the well-posedness of   BSRE is studied. 

\section{Preliminaries}
Let $(\varOmega,\,\mathcal F,\,P)$ be a complete probability  space with  filtration $\{\mathcal F_t\}_{t\geq 0}$, and let   $W$ be an  $m$-dimensional standard  Brownian motion defined on $(\varOmega,\,\mathcal F,\,P)$. We assume ${\mathcal F_t}=\sigma(W_s, 0\le s\le t)$. 
Let $H\in(1/2,\,1)$ be a fixed number and denote 
\begin{equation*}
Z_H(t,\,s)=\kappa_H\left[\left(\frac{t}{s}\right)^{H-\frac{1}{2}}(t-s)^{H-\frac{1}{2}}-\left(H-\frac{1}{2}\right)s^{\frac{1}{2}-H}\int_s^t u^{H-\frac{3}{2}}(u-s)^{H-\frac{1}{2}}du\right]
\end{equation*}
with
\begin{equation*}
\kappa_H=\sqrt{\frac{2H\varGamma(\frac{3}{2}-H)}{\varGamma(H+\frac{1}{2})\varGamma(2-2H)}}.
\end{equation*}
Define 
\begin{align}
W_j^H(t)=\int_0^t Z_H(t,\,s) dW_j(s), \quad 1\le j\le m,
\end{align}
then $W_t^H=\left(W_1^H(t),\cdots,W_m^H(t) , t\in[0,T]\right)$  
is an   $m$-dimensional fractional Brownian motion %$W_t^H=\left(W_1^H(t),\cdots,W_m^H(t)\right), t\in[0,T]$  
of Hurst parameter $H$ (e.g. \cite{hu2005integral}).  Namely, it is a continuous, mean 0 Gaussian process with   covariance
$$ \mathbb{E}[W_i^H(t)W_j^H(s)]=\frac{1}{2}\delta_{ij}(t^{2H}+s^{2H}-| t-s|^{2H})\,,\quad  \forall \ s, t \ge 0\,.   $$

%Let a fractional Brownian motion $B_t^H=\left(B_1^H(t),\cdots,B_m^H(t)\right), t\in[0,T]$ and a standard Brownian motion $B_t=\left(B_1(t),\cdots,B_m(t)\right), t\in[0,T]$ satisfy
%\begin{align}
%B_j^H(t)=\int_0^t Z_H(t,\,s) dB_j(s), \quad 1\le j\le m,
%\end{align}
%where 

We call $W_t$   the underlying Brownian motion of $W_t^H$.
Let the operators  $\mathbf{\Gamma}_{H,T}^*$ and $\mathbf{B}_{H,T}^*$ (\cite{hu2005integral}) be defined by 
\begin{align}
\mathbf{\Gamma}_{H,T}^*f(t)=\left(H-\frac{1}{2}\right)\kappa_H t^{\frac{1}{2}-H}\int_t^T u^{H-\frac{1}{2}}(u-t)^{H-\frac{3}{2}}f(u)du, \quad 0\leq t\leq T
\label{2.3}
\end{align}
and
\begin{align*}
\mathbf{B}_{H,T}^*f(t)=-\frac{2H\kappa_1}{\kappa_H}t^{\frac{1}{2}-H}\frac{d}{dt}\int_t^T (u-t)^{\frac{1}{2}-H}u^{H-\frac{1}{2}}f(u)du,
\end{align*}
where $f$ is a deterministic function and
\begin{align*}
\kappa_1=\frac{1}{2H\Gamma(H-\frac{1}{2})\Gamma(\frac{3}{2}-H)}.
\end{align*}
Furthermore, we have (Definition 6.1 in \cite{hu2005integral})
\begin{align}
\int_0^T f(t)dW_j^H(t):=\int_0^T (\mathbf{\Gamma}_{H,T}^*f)(t)dW_j(t)
\label{2.4}
\end{align}
and
\begin{align}
\int_0^T f(t)dW_j(t):=\int_0^T (\mathbf{B}_{H,T}^*f)(t)dW^H_j(t).
\label{e.2.4} 
\end{align}

It is obvious   that the filtration generated by the fractional Brownian motion $\{W_t^H\}_{t\in[0,\,T]}$ coincides with the one generated by the underlying Brownian motion $\{W_t\}_{t\in[0,\,T]}$, because if we take $f=1$, then \eqref{e.2.4} gives 
\begin{align*}
W_j (s)=\int_0^s (\mathbf{B}_{H,T}^*1)(t)dW^H_j(t) \,, 
\end{align*}
which  
is $\mathcal{F}_s^{W^H}$-measurable for $\forall\, 0\le s \le T$. 
%In the same way, by using operator $\mathbb{B}^*_{H,T}$, we show that
%the filtration generated by $\{W_t\}_{t\in[0,\,T]}$ coincides with the one generated by $\{W_t^H\}_{t\in[0,\,T]}$. So
%\begin{align*}
%\mathcal{F}_t=\sigma(B^H_s,\ 0\leq s\leq t)=\sigma(W_s,\ 0\leq s\leq t).
%\end{align*}

Next, we introduce the Malliavin derivatives \cite{hu2005integral}. Let $\xi_1,\cdots,\xi_k,\cdots$ be an orthonormal basis of $L^2([0,\,T])$ such that $\xi_k$, $k=1,2,\cdots$ are smooth functions on $[0,\,T]$.   Define $\eta_j=\mathbb{B}_{H,T}\xi_j$ and denote by $\mathcal{P}^H$ the set of all polynomial functional $\tilde{\eta}_{j,l}=\int_0^T\eta_j(t)dW_l^H(t)$. The elements G in $\mathcal{P}^H$ have the following form
\begin{align*}
G(\omega)=g(\tilde{\eta}_{j_1,l_1},...,\tilde{\eta}_{j_n,l_n}),
\end{align*}
where $g$ is a polynomial of $n$ variables. We define its Malliavin derivative $D_s^{H,l}G$ by
\begin{equation*}
D_s^{H,l}G=\sum_{i=1}^{n}\frac{\partial g}{\partial x_i}\left(\tilde{\eta}_{j_1,l_1},...,\tilde{\eta}_{j_n,l_n}\right)\eta_{j_i}(s)I_{\{j_i=l\}},\quad \forall\  0\le s\le T.
\end{equation*}

Another useful Malliavin derivative $\mathbb{D}_t^{H,j}G$ is defined as follows 
\begin{equation*}
\mathbb{D}_t^{H,j} G=\int_0^T \phi(s,\,t)D_s^{H,j} Gds, \ \ 0\leq t\leq T,
\end{equation*}
where 
\begin{align*}
\phi(s,\,t)=H(2H-1)\lvert s-t\lvert^{2H-2}.
\end{align*}
The following equality will be used: 
\begin{align}\label{2.4}
\int_a^bX_s\circ dW_s^H=\int_a^bX_sdW_s^H+\sum_{j=1}^m\int_a^b\mathbb{D}_s^{H,j}X_s ds.
\end{align}

\section{Dynamic programming principle}
In this section, we consider the following control problem: the state equation is 
\begin{align}\label{state}
\left\{\begin{array}{l}
 d X_s^{t,x,u}=b\left(s,X_s^{t,x,u},u_s\right)dt+\sigma\left(s,X_s^{t,x,u}\right)\circ dW^H_s\\
\qquad\qquad+\gamma\left(s,X_s^{t,x,u},u_s\right)dW_s,\quad s\in[t,T],\\
 X_t^{t,x,u}=x,
\end{array}\right.
\end{align}
with the cost functional
\begin{align*}
J\left(u;t,x\right)=\mathbb{E}^{\mathcal{F}_t}\left[\int_t^Tf(s,X_s^{t,x,u},u_s)ds+g(X_T^{t,x,u})\right].
\end{align*}
Here $x\in\mathbb{R}^d$,  $W^H$ is an   $m$-dimensional fractional Brownian motion with Hurst parament $H>1/2$, and $W=(W^1,...,W^m)$ is the ``underlying" standard Brownian motion for  $W^H$,  $b,\sigma,\gamma,f,g$ are measurable functions  
\[
\left\{\begin{aligned}
b&: \Omega\times[0,T]\times\mathbb{R}^d\times\mathbb{R}^k\to\mathbb{R}^d;\\
\sigma&:\Omega\times[0,T]\times\mathbb{R}^d\to\mathbb{R}^{d\times m};\\
\gamma&:\Omega\times[0,T]\times\mathbb{R}^d\times\mathbb{R}^k\to\mathbb{R}^{d\times m};\\
f&:\Omega\times[0,T]\times\mathbb{R}^d\times\mathbb{R}^k\to\mathbb{R};\\
g&:\Omega\times\mathbb{R}^d\to\mathbb{R}.
\end{aligned} \right. 
\]
They satisfy  the following assumptions:
\begin{enumerate}
	\item[(H3.1)]  $b(\cdot,x,u),\sigma(\cdot,x),\gamma(\cdot,x,u),f(\cdot,x,u)$ are adapted processes and are in $L^2_\mathcal{F}(0,T)$ for all $(x,u)\in\mathbb{R}^d\times\mathbb{R}^k$.

\item[(H3.2)] $b,\sigma,\gamma$ are Lipschitz continuous w.r.t. $(x,u)$. This means that there exists some $L>0$, such that
\begin{align*}
&|b(t,x_1,u_1)-b(t,x_2,u_2)|+\|\sigma(t,x_1)-\sigma(t,x_2)\|\\
&\qquad +\|\gamma(t,x_1,u_1)-\gamma(t,x_2,u_2)\|\le L\left(|x_1-x_2|+|u_1-u_2|\right),\quad a.s.
\end{align*}
for all $t\in[0,T],\, x_1,x_2\in\mathbb{R}^d, \,u_1,u_2\in\mathbb{R}^k$.

\item[(H3.3)] $\sigma(t,x)$ is H\"{o}lder continuous in time uniformly in $x\in\mathbb{R}^d$. There exist  some constants $1-H<\gamma<1$, and $L>0$, such that
\begin{align*}
\|\sigma(t,x)-\sigma(s,x)\|\le L|t-s|^\gamma,\quad
\forall t,s\in[0,T],\quad \forall x\in\mathbb{R}^d.
\end{align*}

\item[(H3.4)] For $f,g$, there exists some $L>0$, such that
\begin{align*}
&|f(t,x_1,u_1)-f(t,x_2,u_2)|+|g(x_1)-g(x_2)|\\
&\qquad\quad \le L(1+|x_1|+|x_2|)|x_1-x_2|+L(1+|u_1|+|u_2|)|u_1-u_2|, \quad a.s.
\end{align*}
for all $x_1,x_2\in\mathbb{R}^d,\, u_1,u_2\in\mathbb{R}^k$.
\end{enumerate}

Let  $\mathcal{U}[a,b]$ denote the set of admissible control process
$u=(u_t)_{a\le t\le b}$ which is ${\mathcal F}_t$-adapted, $\beta$-H\"older continuous for some $\beta>1/2$, 
and takes value  in a given domain $\textbf{D}\subseteq \mathbb{R}^k$ satisfying $\mathbb{E}\int_a^b |u_t|^2 dt <+\infty$.
%\footnote{Need to explain the controlled system has a solution? }
 
Under the above assumptions, the existence and uniqueness  of the controlled system (\ref{state})  are known (e.g.  Silva et al. \cite{da2018mixed}). 

Define the value function by
\begin{align*}
V(t,x)=\inf_{u\in\mathcal{U}[t,T]}J(u;t,x),\quad (t,x)\in[0,T]\times\mathbb{R}^d.
\end{align*}

Now we discuss the DPP for our stochastic control problem. Similar to the classical case, the following property holds for the value function.

\begin{proposition}
Assume that (H3.1)-(H3.4) hold, then the value function obeys  the following DPP: for all $0\le t\le t+\delta\le T$ and $x\in\mathbb{R}^d$,
\begin{align*}
V(t,x)=\inf_{u\in\mathcal{U}[t,t+\delta]}\mathbb{E}^{\mathcal{F}_t}\left[\int_t^{t+\delta}f(s,X_s^{t,x,u},u_s)ds+V\left(t+\delta,X_{t+\delta}^{t,x,u}\right)\right].
\end{align*}
\end{proposition}
\begin{proof}
The idea of the proof is classical. Through the definition of the value function, we have
\begin{align*}\label{vtx}
V(t,x)=\inf_{u\in\mathcal{U}[t,T]}\mathbb{E}^{\mathcal{F}_t}\left[\int_t^{t+\delta}f(s,X_s^{t,x,u},u_s)ds+J\left(u;t+\delta,X_{t+\delta}^{t,x,u}\right)\right].
\end{align*}
Then it is obvious  that
\begin{align} 
V(t,x)&\ge\inf_{u\in\mathcal{U}[t,T]}\mathbb{E}^{\mathcal{F}_t}\left[\int_t^{t+\delta}f(s,X_s^{t,x,u},u_s)ds+V\left(t+\delta,X_{t+\delta}^{t,x,u}\right)\right]\notag \\
&=\inf_{u\in\mathcal{U}[t,t+\delta]}\mathbb{E}^{\mathcal{F}_t}\left[\int_t^{t+\delta}f(s,X_s^{t,x,u},u_s)ds
+V\left(t+\delta,X_{t+\delta}^{t,x,u}\right)\right].\label{vtxge}
\end{align}
On the other hand, for all $(t,x)\in[0,T]\times\mathbb{R}^d$ and $\varepsilon>0$, there exists a $\tilde{u}^{t,x}\in\mathcal{U}[t,T]$, such that
\begin{align*}
J(\tilde{u}^{t,x};t,x)\le \inf_{u\in\mathcal{U}[t,T]}J(u;t,x)+\varepsilon=V(t,x)+\varepsilon, \, a.s.
\end{align*}
Thus,  
\begin{align}\label{vtxle}
V(t,x)&= \inf_{u\in\mathcal{U}[t,T]}\mathbb{E}^{\mathcal{F}_t}\left[\int_t^{t+\delta}f(s,X_s^{t,x,u},u_s)ds+J\left(\tilde{u}^{t+\delta,X_{t+\delta}^{t,x,u}};t+\delta,X_{t+\delta}^{t,x,u}\right)\right]\\
&\le\inf_{u\in\mathcal{U}[t,t+\delta]}\mathbb{E}^{\mathcal{F}_t}\left[\int_t^{t+\delta}f(s,X_s^{t,x,u},u_s)ds+V\left(t+\delta,X_{t+\delta}^{t,x,u}\right)\right]+\varepsilon.\notag
\end{align}
By the arbitrariness  of $\varepsilon$, and the  equations (\ref{vtxge}) and (\ref{vtxle}), we complete the proof.
\end{proof}

Next,  we give the main result of this section.
\begin{theorem}\label{dpp}
Assume that the value function has the following representation: 
%\footnote{Is this decomposition unique? Namely, Are 
 %$\Phi(t,\cdot), \Psi(t,\cdot), Q^i(t,\cdot)$ uniquely determined by %$V(t,x)$?}
\begin{align}
V(t,x)=g(x)+\int_t^T\Phi(s,x)ds+\int_t^T\Psi(s,x)\circ dW^H_s-\sum_{i=1}^m\int_t^TQ^i(s,x)dW^i_s,\label{e.3.4} 
\end{align}
for $t\in[0,T]$. Here $\Phi(t,\cdot), \Psi(t,\cdot), Q^i(t,\cdot)$ are $\mathcal{F}_t\times\mathcal{B}(\mathbb{R}^d)$ measurable functions with  values in $\mathbb{R}, \mathbb{R}^{1\times m},\mathbb{R}$, respectively, for each $t\in[0,T]$. Moreover, we make   the following assumptions. 
\begin{itemize}
    \item[(1)]  $(t,x)\longmapsto V(t,x)$ is jointly continuous a.s.,
    \item[(2)] $x\longmapsto V(t,x)$ is $C^2$ for each $t\in[0,T]$ a.s.,
    \item[(3)] $x\longmapsto \Phi(t,x),\Psi(t,x)$ are continuous for each $t\in[0,T]$ a.s.,
    \item[(4)] $x\longmapsto Q^i(t,x)$ are $C^1$ for each $t\in[0,T]$ a.s.,
    \item[(5)]  There exits $L\in L_\mathcal{F}^2(0,T;\mathbb{R}^+)$, such that
    \end{itemize}
\begin{align*}
|V(t, x)|,|g(x)|,&|\Phi(t, x)|,|\Psi(t, x)|,|Q^i(t,x)| \le L_t\left(1+|x|^2\right), \\
&\left| V_x(t, x)\right|,\left\| Q^i_x(t, x)\right\| \leq L_t(1+|x|), \\
&\quad\left\|V_{xx}(t, x)\right\| \leq L_t, \\
|\Phi(t, x)-\Phi(t,y)|&+|\Psi(t,x)-\Psi(t,y)| \leq L_t(1+|x|+|y|)|x-y|.
\end{align*}
Then the value function $V$, together with $Q$, is a solution pair to  the following stochastic HJB equation:
\begin{equation}\label{hjb}
    \left\{\begin{array}{l}
 -dV(t,x)=\inf_uG(t,x,V_x(t,x),V_{xx}(t,x),Q_x(t,x),u)dt\\
\qquad\qquad\qquad+\langle \sigma(t,x),V_x(t,x)\rangle\circ dW^H_t-\sum_{i=1}^mQ^i(t,x)dW^i_t,\\
 V(T,x)=g(x),
\end{array}\right.
\end{equation}
where $Q_x=(Q^1_x,...,Q_x^m)$, and
\begin{align*}
G(t,x,p,P,A,u)=&\langle b(t,x,u),p\rangle+f(t,x,u)+\sum_{i=1}^m\langle\gamma^i(t,x,u),A^i\rangle\\
&+\frac{1}{2} \hbox{\rm tr} \langle P\gamma(t,x,u),\gamma(t,x,u)\rangle.
\end{align*}
\end{theorem}
\begin{proof}  Applying the It$\hat{\rm o}$-Ventzell 
 formula (see in \cite{ad1965equations,ustunel1982some,kunita2006stochastic})
% \footnote{give a reference} 
to $V(t,X_t)$, with 
% \footnote{why $\Psi(t,x)=\langle \sigma(t,x),V_x(t,x)\rangle$?
 %You need the decomposition of the form \eqref{e.3.4} is unique?}  
  $\Psi(t,x)=\langle \sigma(t,x),V_x(t,x)\rangle$, we have
\begin{align*}
&dV(t,X_t)\\=&-\Phi(t,X_t)dt-\Psi(t,X_t)\circ dW_t^H+\sum_{i=1}^mQ^i(t,X_t)dW^i_t\\
&+\langle b(t,X_t,u_t),V_x(t,X_t)\rangle dt+\langle\sigma(t,X_t),V_x(t,X_t)\rangle\circ dW_t^H\\
&+\sum_{i=1}^m\langle\gamma^i(t,X_t,u_t),Q^i_x(t,X_t)\rangle dt+\frac{1}{2}\tr\langle V_{xx}(t,X_t)\gamma(t,X_t,u_t),\gamma(t,X_t,u_t)\rangle dt\\
&+\langle V_x(t,X_t),\gamma(t,X_t,u_t)\rangle dW_t\\
=&[-\Phi(t,X_t)+G(t,X_t,V_x(t,X_t),V_{xx}(t,X_t),Q_x(t,X_t),u_t)-f(t,X_t,u_t)]dt\\
&+[Q(t,X_t)+V_x(t,X_t)\gamma(t,X_t,u_t)]dW_t.
\end{align*}
So, 
\begin{align*}
V(t,x)=&\mathbb{E}^{\mathcal{F}_t}\left[V(t+\delta,X_{t+\delta}^{t,x,u})+\int_t^{t+\delta}f(s,X_s^{t,x,u},u_s)ds\right]\\
&+\mathbb{E}^{\mathcal{F}_t}\int_{t}^{t+\delta}\Delta(s,X_s^{t,x,u},u_s)ds,
\end{align*}
where
\begin{align}\label{Delta}
\Delta(s,x,u)=\Phi(s,x)-G(s,x,V_x(s,x),V_{xx}(s,x),Q_x(s,x),u).
\end{align}
From the DPP, we have
\begin{align*}
V(t,x)\le \mathbb{E}^{\mathcal{F}_t}\left[V(t+\delta,X_{t+\delta}^{t,x,u})+\int_t^{t+\delta}f(s,X_s^{t,x,u},u_s)ds\right],
\end{align*}
which implies
\begin{align*}
\mathbb{E}^{\mathcal{F}_t}\int_{t}^{t+\delta}\Delta(s,X_s^{t,x,u},u_s)ds\le 0.
\end{align*}
Due to the arbitrariness of $t,\delta,x$, we have 
\begin{align*}
\Delta(t,x,u_t)\le 0, \quad\forall(t,x)\in[0,T]\times\mathbb{R}^d, \, a.s.
\end{align*}
From  the arbitrariness of $u$, it follows
\begin{align}\label{3.5}
\sup_{u}\Delta(t,x,u)\le 0, \quad\forall(t,x)\in[0,T]\times\mathbb{R}^d, \, a.s.
\end{align}
On the other hand, for all $t,\delta,x$ and $\varepsilon>0$, there exists $\tilde{u}\in\mathcal{U}[t,t+\delta]$, such that
\begin{align*}
\mathbb{E}^{\mathcal{F}_t}\left[V(t+\delta,X_{t+\delta}^{t,x,\tilde{u}})+\int_t^{t+\delta}f(s,X_s^{t,x,\tilde{u}},\tilde{u}_s)ds\right]\le V(t,x)+\varepsilon,
\end{align*}
which means
\begin{align*}
\mathbb{E}^{\mathcal{F}_t}\int_{t}^{t+\delta}\Delta(s,X_s^{t,x,\tilde{u}},\tilde{u}_s)ds\ge -\varepsilon.
\end{align*}
Due to the arbitrariness of $t,\delta,x$ and $\varepsilon$, combining with equation (\ref{3.5}), we have
\begin{align*}
\sup_{u}\Delta(t,x,u)= 0, \quad\forall(t,x)\in[0,T]\times\mathbb{R}^d, \, a.s.
\end{align*}
This gives 
\begin{align*}
\Phi(t,x)=\inf_uG(t,x,V_x(t,x),V_{xx}(t,x),Q_x(t,x),u) \,, 
\end{align*}
and hence the proof is completed. 
\end{proof}

It remains to show the uniqueness of the decomposition  (\ref{e.3.4}). We have the following theorem.
\begin{theorem} 
Let $V(t), \Phi(t),\Psi(t),Q(t)$ be adapted  random processes   
satisfying 
\begin{align}
V(t)=V(T) +\int_t^T\Phi(s )ds+\int_t^T\Psi(s )\circ dW^H_s-\sum_{i=1}^m\int_t^TQ^i(s ) dW^i_s,\,\forall t\in[0,T].
\label{e.3.4a} 
\end{align}

   Then $\Phi,\Psi,Q$ are uniquely determined by $V$.
\end{theorem}
\begin{proof}
Let both $(\Phi,\Psi,Q)$ and $(\tilde{\Phi},\tilde{\Psi},\tilde{Q})$ satisfy the equation (\ref{e.3.4a}) with the same $V$. Denote $\Delta\Phi=\tilde{\Phi}-\Phi, \Delta\Psi=\tilde{\Psi}-\Psi, \Delta Q=\tilde{ Q}- Q$.   Then we have
\begin{align}\label{deltaV}
0=\Delta V(t):=\int_t^T\Delta\Phi(s )ds+\sum_{i=1}^m\int_t^T\Delta\Psi^i(s )\circ dW^{H,i}_s-\sum_{i=1}^m\int_t^T\Delta Q^i(s )dW^i_s.
\end{align}
Applying  It$\hat{\rm o}$'s formula to $\left|\Delta V(t )\right|^2$, we have
\begin{align*}
\mathbb{E}\sum_{i=1}^m\int_0^T\left|\Delta Q^{i}(t )\right|^2dt=&\mathbb{E}\left|\Delta V(T )\right|^2-\mathbb{E}\left|\Delta V(0 )\right|^2\\
&+\mathbb{E}\int_0^T\Delta V(t )\Delta \Phi( )dt\\
&+\mathbb{E}\sum_{i=1}^m\int_0^T\Delta V(t  )\Delta\Psi(t )\circ dW_t^{H,i}\\
=&0,
\end{align*}
because $\Delta V(t)=0$,  which shows the uniqueness of $Q$. Through equation (\ref{deltaV}), it is clear that
\begin{align*}
\mathbb{E}^{\mathcal{F}_a}\left[\int_a^b\Delta\Phi(s )ds+\sum_{i=1}^m\int_a^b\Delta\Psi^i(s )\circ dW_s^{H,i}\right]=0,\quad a.s., \quad \forall\  0\le a<b\le T.
\end{align*}
According to Lemma 5.6 of \cite{han2013maximum}, we have
\begin{align*}
\Delta\Phi(a )=\Delta\Psi(a )=0,\quad a.s.,\quad \forall a\in[0,T],
\end{align*}
since $\Delta\Phi(a )$ and $\Delta\Psi(a )$ are $\mathcal{F}_a$-measurable. The proof is complete.
\end{proof}

We have proved that the value function satisfies HJB equation (\ref{hjb}),  which will imply  the following  stochastic verification theorem.
\begin{theorem}\label{verification}
Let    $(\mathbb{V,Q})$ be the solution to equation (\ref{hjb}) and satisfy the assumptions in Theorem \ref{dpp}. Then for any $(t,x)\in[0,T]\times\mathbb{R}^d$ and $u\in\mathcal{U}[t,T]$, we have
\begin{align*}
\mathbb{V}(t,x)\le J(u;t,x).
\end{align*}
Moreover, if there exists a $u^*\in\mathcal{U}$, such that 
\begin{align*}
&G(s,X_s^{t,x,u^*},\mathbb{V}_x(s,X_s^{t,x,u^*}),\mathbb{V}_{xx}(s,X_s^{t,x,u^*}),\mathbb{Q}_x(s,X_s^{t,x,u^*}),u_s^*)\\
=&\inf_{v}G(s,X_s^{t,x,u^*},\mathbb{V}_x(s,X_s^{t,x,u^*}),\mathbb{V}_{xx}(s,X_s^{t,x,u^*}),\mathbb{Q}_x(s,X_s^{t,x,u^*}),v),
\end{align*}
a.s., for all $t\le s\le T$, $x\in\mathbb{R}^d$. Then $u^*$ is the optimal control. 
\end{theorem}
\begin{proof}
Applying  the It$\hat{\rm o}$ formula to $\mathbb{V}(s,X_s^{t,x,u})$ and taking conditional expectation yield 
\begin{align*}
\mathbb{V}(t,x)=&\mathbb{E}^{\mathcal{F}_t}\left[\int_t^Tf(s,X_s^{t,x,u})ds+g(X_T^{t,x,u})\right]\\
&-\mathbb{E}^{\mathcal{\mathcal{F}}_t}\int_t^TG(s,X_s^{t,x,u},\mathbb{V}_x(s,X_s^{t,x,u}),\mathbb{V}_{xx}(s,X_s^{t,x,u}),\mathbb{Q}_x(s,X_s^{t,x,u}),u_s)ds\\
&+\mathbb{E}^{\mathcal{\mathcal{F}}_t}\int_t^T\inf_vG(s,X_s^{t,x,u},\mathbb{V}_x(s,X_s^{t,x,u}),\mathbb{V}_{xx}(s,X_s^{t,x,u}),\mathbb{Q}_x(s,X_s^{t,x,u}),v)ds\\
\le&J(u;t,x),\quad a.s., \quad \forall(t,x)\in[0,T]\times\mathbb{R}^d,\, u\in\mathcal{U}[t,T].
\end{align*}

In the same way, we have
\begin{align*}
\mathbb{V}(t,x)=\mathbb{E}^{\mathcal{F}_t}\left[\int_t^Tf(s,X_s^{t,x,u^*})ds+g(X_T^{t,x,u^*})\right]=J(u^*;t,x)\le\inf_{u\in\mathcal{U}[t,T]}J(u;t,x),
\end{align*}
which means $u^*$ is the optimal control.
\end{proof}

Consider the case with no control and $f=0$. Then the verification theorem is reduced to the following Feynman-Kac formula.

\begin{theorem}
Let $X$ be the solution to the following SDE: 
\begin{align*}
\left\{\begin{array}{l}
 d X_s^{t,x}=b\left(s,X_s^{t,x}\right)dt+\sigma\left(s,X_s^{t,x}\right)\circ dW^H_s\\
\qquad\qquad+\gamma\left(s,X_s^{t,x}\right)dW_s,\quad s\in[t,T],\\
 X_t^{t,x}=x,
\end{array}\right.
\end{align*}
where the coefficients satisfy regular assumptions similar to the control problem. Then 
\begin{align*}
\tilde{V}(t,x):=\mathbb{E}^{\mathcal{F}_t}g(X_T^{t,x})
\end{align*}
satisfies the following backward stochastic partial differential equations:
\begin{equation}\label{Feymankac}
    \left\{\begin{array}{l}
 -d\tilde{V}(t,x)=\left[\langle b(t,x,u),\tilde{V}_x(t,x)\rangle+\sum_{i=1}^m\langle\gamma^i(t,x),\tilde{Q}_x^i(t,x)\rangle\right] dt\\
 \qquad\qquad\qquad+\frac{1}{2}tr\langle \tilde{V}_{xx}(t,x)\gamma(t,x),\gamma(t,x)\rangle dt\\
\qquad\qquad\qquad+\langle \sigma(t,x),\tilde{V}_x(t,x)\rangle\circ dW^H_t-\sum_{i=1}^m\tilde{Q}^i(t,x)dW^i_t,\\
 \tilde{V}(T,x)=g(x),
\end{array}\right.
\end{equation}
\end{theorem}

\begin{example}
Consider the  $1$-dimensional case. Let $b=\gamma=0$, $\sigma=1$ and $g(x)=x$. Then
\begin{align*}
X_s^{t,x}=x+W_s^H-W_t^H,\quad s\in[0,T],
\end{align*}
and
\begin{align*}
\tilde{V}(t,x)=\mathbb{E}^{\mathcal{F}_t}X_T^{t,x}=x+\mathbb{E}^{\mathcal{F}_t}W_T^H-W_t^H.
\end{align*}
Notice that $\mathbb{E}^{\mathcal{F}_t}W_T^H$ is a square integrable martingale. Then there exists a unique $Q\in L_\mathcal{F}^2(0,T;\mathbb{R})$, such that 
\begin{align*}
\mathbb{E}^{\mathcal{F}_t}W_T^H=\int_0^tQ(s)dW_s,\quad \forall t\in[0,T].
\end{align*}
Then we rewrite $\tilde{V}$ as
\begin{align*}
\tilde{V}(t,x)=x-\int_0^t1\circ dW_s^H+\int_0^tQ(s)dW_s.
\end{align*}
So that
\begin{align*}
\tilde{V}(t,x)=\tilde{V}(T,x)+\int_t^T1\circ dW_s^H-\int_t^TQ(s)dW_s,
\end{align*}
which implies BSPDE (\ref{Feymankac}) with the terminal condition $\tilde{V}(T,x)=x$.
\end{example}

To end  this section, we introduce the relationship between the DPP and stochastic maximum principle for our control problem. 
\begin{theorem}
Assume that the optimal control $u^*$ exists. Let $(p,q)$ be the solution to the following adjoint BSDE
\begin{equation}\label{adjoint}
    \left\{\begin{array}{l}
 -dp_s^{t,x}=[W_x^*(s)p_s^{t,x}+\gamma^*_x(s)q_s^{t,x}+f_x^*(s)]ds\\
 \qquad\qquad\qquad+\sigma_x^*(s)p_s^{t,x}\circ dW_s^H-q_s^{t,x}dW_s,\\
p_T^{t,x}=g_x(X_T^{t,x,u^*}),
\end{array}\right.
\end{equation}
where we use the notation 
\begin{align*}
l^*(s)=l(s,X_s^{t,x,u^*},u^*_s),
\end{align*}
for $l=b_x,\sigma_x,\gamma_x,f_x$.

Then the following equality hold:
\begin{align*}
p_s^{t,x}=V_x(s,X_s^{t,x,u^*}),\quad q_s^{t,x}=Q_x(s,X_s^{t,x,u^*}).
\end{align*}
\end{theorem}
\begin{proof}
Applying the  It$\hat{\rm o}$-Ventzell formula to $V_x(s,X_s^{t,x,u^*})$, and combining the existence and uniqueness result of BSDE (\ref{adjoint}) (We refer to Han and Sun \cite{yuecai2018solutions}), one   completes  the proof.
\end{proof}
\section{Linear quadratic optimal control} In this section, we apply  the  DDP obtained in previous  section to the LQ case. In one   state space dimension,  
and Markov control case,  the LQ control problem was solved in 
in \cite{hu2005stochastic} and \cite{duncan2013linear} in the absence of  Brownian motion ($\gamma=0$).  Our DDP method can be applied to multi space dimensions. But for notational  simplexity  without losing generality, we assume $d=k=m=1$. The state equation is 
\begin{align*}
\left\{\begin{array}{l}
 d X_s^{t,x,u}=\left(A_sX_s^{t,x,u}+B_su_s\right) dt+(C_sX_s^{t,x,u}+D_su_s)dW_s\\
 \qquad\qquad\qquad+F_sX_s^{t,x,u}\circ dW^H_s,\quad s\in[t,T],\\
 X_t^{t,x,u}=x\,. 
\end{array}\right.
\end{align*}
We want to  minimize the cost functional
\begin{align*}
J(u;t,x)=\mathbb{E}^{\mathcal{F}_t}\left[\int_t^T\left(S_s|X_s^{t,x,u}|^2+R_su_s^2 \right)ds+G|X_T^{t,x,u}|^2\right].
\end{align*}
We make  following assumptions on  the coefficients:
\begin{enumerate}
	\item[(H4.1)]  $A,B,C,D,F,S,R$ are adapted processes, such that $S_s\ge 0,\, R_s>0,\,a.s.$ for all $s\in[0,T]$. $G$ is an  $\mathcal{F}_T$-measurable random variable, such that $G\ge 0$, a.s.

\item[(H4.2)] Let $\tilde{\mathcal{Y}}_p$   defined by  
\begin{align*}
\tilde{\mathcal{Y}}_p=\sup_{0\le t\le s\le T}\exp\left\{p\cdot\left[\int_t^s(A_r-\frac{1}{2}C_r^2)dr+\int_t^sC_rdW_r+\int_t^sF_r\circ dB_r^H\right]\right\} 
\end{align*}
exists and  $\mathbb{E}\tilde{\mathcal{Y}}_p<+\infty$ for all $p\in[-8,8]$.

\item[(H4.3)]
Let $\tilde{\mathcal{M}}$  defined by  
\begin{align*}
\tilde{\mathcal{M}}=\tilde{\mathcal{Y}}_4\times\left(\int_0^TS_sds+G\right)^2 
\end{align*}
satisfies  $\mathbb{E}\tilde{\mathcal{M}}<+\infty$, and
\begin{align*}
\mathbb{E}\left[\tilde{\mathcal{Y}}_8\times\left(G^4+\int_0^T\left[\tilde{\mathcal{M}}^2(A_s^2+C_s^2+C_s^4)+\tilde{\mathcal{M}}S_s^2\right]ds\right)\right]<+\infty.
\end{align*}

\item[(H4.4)] $F_t$ is H\"{o}lder continuous. There exists some constants $1-H<\gamma<1$, and $L>0$, such that
\begin{align*}
|F_t-F_s|\le L|t-s|^\gamma,\quad a.s.,\quad
\forall t,s\in[0,T].
\end{align*}
\end{enumerate}

The HJB equation for this  LQ control problem becomes 
\begin{proposition}
Let assumptions (H4.1)-(H4.4) hold. Then the value function satisfies the following HJB equation:
\begin{equation}\label{hjblq}
    \left\{\begin{array}{l}
 -dV(t,x)=-\frac{1}{2}(2R_t+D_t^2V_{xx}(t,x))^{-1}\\\qquad\qquad\qquad\times\left[B_tV_x(t,x)+D_tQ_x(t,x)+C_tD_txV_{xx}(t,x)\right]^2dt\\
 \qquad\qquad\qquad+\left[(A_tV_x(t,x)+C_tQ_x(t,x))x+(S_t+\frac{1}{2}C_t^2V_{xx}(t,x))x^2\right]dt\\
\qquad\qquad\qquad+F_txV_x(t,x)\circ dW^H_t-Q(t,x)dW_t,\\
 V(T,x)=Gx^2.
\end{array}\right.
\end{equation}
Moreover, the optimal control $u^*$ is given by
\begin{align}\label{u*1}
u^*_s
=&-\big(2R_s+D_s^2V_{xx}(s,X_s^{t,x,u^*})\big)^{-1}\\&\qquad\times\left(B_sV_x(s,X_s^{t,x,u^*})+D_sQ_x(s,X_s^{t,x,u^*})+C_sD_sxV_{xx}(s,X_s^{t,x,u^*})\right)\notag.
\end{align}
\end{proposition}

\begin{proof}
We apply  Theorem \ref{dpp} with   $G$  being given   by 
\begin{align*}
G(t,x,p,P,A,u)=&(A_tx+B_tu)p+(C_tx+D_tu)A\\
&+\frac{1}{2}(C_tx+D_tu)^2P+S_tx^2+R_tu^2.
\end{align*}
Then 
\begin{align}\label{4.3}
\Phi(t,x)=&\inf_{u}G(t,x,V_x(t,x),V_{xx}(t,x),Q_x(t,x),u)\\
=&\inf_u\Big[\big(R_t+\frac{1}{2}D_t^2V_{xx}(t,x)\big) u^2\notag\\
&\qquad +(B_tV_x(t,x)+D_tQ_x(t,x)+C_tD_txV_{xx}(t,x))u\notag\\
&\qquad+(A_tV_x(t,x)+C_tQ_x(t,x))x+(S_t+\frac{1}{2}M_t^2V_{xx}(t,x))x^2\Big]\notag\\
=&-\frac{1}{2}(2R_t+D_t^2V_{xx}(t,x))^{-1}(B_tV_x(t,x)+D_tQ_x(t,x)+C_tD_txV_{xx}(t,x))^2\notag\\
&\qquad+(A_tV_x(t,x)+C_tQ_x(t,x))x+(S_t+\frac{1}{2}M_t^2V_{xx}(t,x))x^2,\notag
\end{align}
and
\begin{align*}
u^*_s=&\arg\inf_uG\left(s,X_s^{t,x,u^*},V_x(s,X_s^{t,x,u^*}),V_{xx}(s,X_s^{t,x,u^*}),Q_x(s,X_s^{t,x,u^*}),u\right)\\
=&-\big(2R_s+D_s^2V_{xx}(s,X_s^{t,x,u^*})\big)^{-1}\\&\qquad\times\left(B_sV_x(s,X_s^{t,x,u^*})+D_sQ_x(s,X_s^{t,x,u^*})+C_tD_sxV_{xx}(s,X_s^{t,x,u^*})\right).
\end{align*}
This completes the proof. 
\end{proof}

As in  Tang \cite{tang2015dynamic}, the value function has the following property.
\begin{proposition}\label{p.4.2} 
The value function 
\begin{align*}
V(t,x)=\inf_{u\in\mathcal{U}[t,T]}\mathbb{E}^{\mathcal{F}_t}\left[\int_t^T\left(S_s|X_s^{t,x,u}|^2+R_su_s^2 \right)ds+G|X_T^{t,x,u}|^2\right]
\end{align*}
satisfies
\begin{align*}
V(t,\xi x)=\xi^2V(t,x),\quad\forall (t,x)\in[0,T]\times\mathbb{R},
\end{align*}
and
\begin{align*}
V(t,x+y)+V(t,x-y)=2V(t,x)+2V(t,y).
\end{align*}
Moreover, there exists an adapted process $K$, such that
\begin{align*}
V(t,x)=K_tx^2.
\end{align*}
\end{proposition}
\begin{proof}
The proof follows Section 3 of Tang \cite{tang2015dynamic}.
\end{proof}

Then we give a upper bound of $K$.

\begin{lemma}\label{lem4.3}
Let Assumptions (H4.1)-(H4.4) hold and let $K_t$ be defined in 
Proposition \ref{p.4.2}.   Then 
\begin{align*}
\mathbb{E}\left[\sup_{t\in[0,T]}K_t^2\right]<+\infty.
\end{align*}
\end{lemma}

\begin{proof}
Define $\mathcal{Y}$ by
\begin{align*}
\left\{\begin{array}{l}
 d \mathcal{Y}_{t,s}=A_s\mathcal{Y}_{t,s}ds+C_s\mathcal{Y}_{t,s}dW_s+F_s\mathcal{Y}_{t,s}\circ dB_s^H,\quad s\in[t,T],\\
\mathcal{Y}_{t,t}=1,
\end{array}\right.
\end{align*}
which has a unique solution
\begin{align}\label{mathcaly}
\mathcal{Y}_{t,s}=exp\left\{\int_t^s(A_r-\frac{1}{2}C_r^2)dr+\int_t^sC_rdW_r+\int_t^sF_r\circ dB_r^H\right\}.
\end{align}
It is easy to check $X_s^{t,x,0}=x\mathcal{Y}_{t,s}$, then
\begin{align*}
\sup_{s\in[t,T]}|X_{s}^{t,x,0}|^2=x^2\sup_{s\in[t,T]}|\mathcal{Y}_{t,s}|^2.
\end{align*}
So that
\begin{align*}
K_tx^2=&V(t,x)\le J(0;t,x)\\
=&\mathbb{E}^{\mathcal{F}_t}\left[\int_t^TS_s|X_s^{t,x,0}|^2ds+G|X_T^{t,x,0}|^2\right]\\
\le&\mathbb{E}^{\mathcal{F}_t}\left[\sup_{s\in[t,T]}|X_s^{t,x,0}|^2\times\left(\int_t^TS_sds+G\right)\right],\quad a.s.
\end{align*}
for all $(t,x)\in[0,T]\times\mathbb{R}$.
So, 
\begin{align*}
\mathbb{E}\left[\sup_{t\in[0,T]}K_t^2\right]\le&\mathbb{E}\left[\sup_{0\le t\le s\le T}|\mathcal{Y}_{t,s}|^4\times\left(\int_0^TS_sds+G\right)^2\right]\\
=&\mathbb{E}\left[\tilde{\mathcal{Y}}_4\times\left(\int_0^TS_sds+G\right)^2\right]\\
<&+\infty\,, 
\end{align*}
giving  desired result. 
\end{proof}

The following theorem gives a characterization   of $K$.

\begin{theorem}
Let $V(t,x)=K_tx^2$. Then $K$, together with $\Theta$, is the solution pair to the following backward stochastic Riccati equation:
\begin{equation}\label{bsre}
    \left\{\begin{array}{l}
 -dK_t=[-(R_t+D_t^2K_t)^{-1}(B_tK_t+D_t\Theta_t+C_tD_tK_t)^2\\
 \quad\qquad\quad+(2A_t+C_t^2)K_t+2C_t\Theta_t+S_t]dt\\
 \qquad\qquad+2F_tK_t\circ dW_t^H-\Theta_tdW_t,\quad t\in[0,T],\\
K_T=G.
\end{array}\right.
\end{equation}
Moreover, the optimal control is given by
\begin{align}\label{u*2}
u^*_s
=&-\big(R_s+D_s^2K_s\big)^{-1}\left((B_s+C_sD_s)K_s+D_s\Theta_s\right)X_s^{t,x,u^*}.
\end{align}
\end{theorem}

\begin{proof}
Substituting  $V(t,x)=K_tx^2$ into HJB equation (\ref{hjblq}) yields 
\begin{align*}
-x^2dK_t=&-\frac{1}{2}(2R_t+2D_t^2K_t)^{-1}\\
&\times\left[2B_tK_tx+D_tQ_x(t,x)+2C_tD_tK_tx\right]^2dt\\
&+\left[(2A_tK_tx+C_tQ_x(t,x))x+(S_t+M_t^2K_t)x^2\right]dt\\
&+2F_tK_tx^2\circ dW^H_t-Q(t,x)dW_t.
\end{align*}
Since the above equality holds for all $x\in\mathbb{R}$, $\frac{Q(t,x)}{x^2}$ should be independent with $x$, which means there exists a $\Theta$, such that $Q(t,x)=\Theta_tx^2$. So,  
\begin{align*}
-x^2dK_t=&-(R_t+D_t^2K_t)^{-1}\\
&\times\left[B_tK_tx+D_t\Theta_t+C_tD_tK_t\right]^2x^2dt\\
&+\left[(2A_tK_tx+2C_t\Theta_tx)x+(S_t+C_t^2K_t)x^2\right]dt\\
&+2F_tK_tx^2\circ dW^H_t-\Theta_tx^2dW_t.
\end{align*}
By the arbitrariness of $x$, and combining with $V(T,x)=K_Tx^2=Gx^2$, we obtain the BSRE (\ref{bsre}). Then substituting  $V(t,x)=K_tx^2,\, Q(t,x)=\Theta_tx^2$ to equation (\ref{u*1}), we obtain the optimal control (\ref{u*2}).
\end{proof}

We now give the verification theorem to LQ problem.

\begin{theorem}\label{the4.5}
Let $(\mathbf{K},\mathbf{\Theta})$ be the solution to BSRE (\ref{bsre}). Let $\mathbf{u}$ be  the process given by
\begin{align}\label{4.7}
\mathbf{u}_s=&-\big(R_s+D_s^2\mathbf{K}_s\big)^{-1}\left((B_s+C_sD_s)\mathbf{K}_s+D_s\mathbf{\Theta}_s\right)X_s^{t,x,\mathbf{u}}.
\end{align}
Then $\mathbf{u}$ is the optimal control, and the value function is given by
\begin{align*}
V(t,x)=J(\mathbf{u};t,x)=\mathbf{K}_tx^2.
\end{align*}
\end{theorem}

\begin{proof}
Let $u\in\mathcal{U}[t,T]$ be any admissible control. We first show $\mathbf{K}_tx^2$ is a lower bound of $J(u;t,x)$ w.r.t. $u$. Through It$\hat{\rm o}$'s formula, we have 
\begin{align*}
&\mathbb{E}^{\mathcal{F}_t}\left[G|X_T^{t,x,u}|^2\right]\\=&\mathbf{K}_tx^2+\mathbb{E}^{\mathcal{F}_t}\int_t^Td\left(\mathbf{K}_s|X_s^{t,x,u}|^2\right)\\
=&\mathbf{K}_tx^2+\mathbb{E}^{\mathcal{F}_t}\int_t^TD_s^2\mathbf{K}_su_s^2ds-\mathbb{E}^{\mathcal{F}_t}\int_t^TS_s\left|X_s^{t,x,u}\right|^2ds\\
&\quad +2\mathbb{E}^{\mathcal{F}_t}\int_t^T\left(B_s\mathbf{K}_s+C_sD_s\mathbf{K}_s+D_s\mathbf{\Theta}_s\right)X_s^{t,x,u}u_sds\\
&\quad +\mathbb{E}^{\mathcal{F}_t}\int_t^T\left(R_s+D_s^2\mathbf{K}_s\right)^{-1}\left((B_s+C_sD_s)K_s+D_s\mathbf{\Theta}_s\right)^2\left|X_s^{t,x,u}\right|^2ds.
\end{align*}
So,   
\begin{align}\label{4.8}
\mathbf{K}_tx^2=&\mathbb{E}^{\mathcal{F}_t}\left[\int_t^T\left(S_s\left|X_s^{t,x,u}\right|^2+R_su_s^2\right)ds+G\left|X_T^{t,x,u}\right|^2\right]\\
&-\mathbb{E}^{\mathcal{F}_t}\int_t^T\left(R_s+D_s^2\mathbf{K}_s\right)u_s^2ds\notag\\
&-2\mathbb{E}^{\mathcal{F}_t}\int_t^T\left(B_s\mathbf{K}_s+C_sD_s\mathbf{K}_s+D_s\mathbf{\Theta}_s\right)X_s^{t,x,u}u_sds\notag\\
&-\mathbb{E}^{\mathcal{F}_t}\int_t^T\left(R_s+D_s^2\mathbf{K}_s\right)^{-1}\left(B_s\mathbf{K}_s+C_sD_s\mathbf{K}_s+D_s\mathbf{\Theta}_s\right)^2\left|X_s^{t,x,u}\right|^2ds\notag\\
=&J(u;t,x)-\mathbb{E}^{\mathcal{F}_t}\int_t^T\left(a_su_s^2+b_sX_s^{t,x,u}u_s+(2a_s)^{-1}b_s^2\left|X_s^{t,x,u}\right|^2\right)ds\notag\\
\le&J(u;t,x),\quad \forall (t,x)\in[0,T]\times\mathbb{R},\quad u\in\mathcal{U}[t,T],\notag
\end{align}
where
\begin{align*}
a_s=R_s+D_s^2\mathbf{K}_s>0,\quad b_s=B_s\mathbf{K}_s+C_sD_s\mathbf{K}_s+D_s\mathbf{\Theta}_s.
\end{align*}
In the same way, letting 
$\mathbf{u}$ be given by equation (\ref{4.7}), we have
\begin{align}\label{4.9}
\mathbf{K}_tx^2=J(\mathbf{u};t,x).
\end{align}
Equations (\ref{4.8}) and (\ref{4.9}) mean 
that $\mathbf{u}$ given by equation (\ref{4.7}) is the optimal control, and 
\begin{align*}
V(t,x)=\inf_{u\in\mathcal{U}[t,T]}J(u;t,x)=J(\mathbf{u};t,x)=\mathbf{K}_tx^2\,, 
\end{align*}
proving the theorem. 
\end{proof}

\section{Solution to backward stochastic Riccati equation} In this section, we study the well-posedness of the solution to BSRE (\ref{bsre}).
We have proved the existence of $K$ in Lemma \ref{lem4.3}, we now give the existence result of $\Theta$.

\begin{theorem}\label{the5.1}
Let assumptions (H4.1)-(H4.4) hold. Then $\Theta\in L_{\mathcal{F}}^2(0,T;\mathbb{R})$.
\end{theorem}

\begin{proof}
Recall the notation $X=X^{0,x,0}$. Using It$\hat{\rm o}$'s formula, we have
\begin{align*}
d\left(K_tX_t^2\right)=&dV(t,X_t)\\
=&\left[-\tilde{\Phi}(t,X_t)+(2A_t+C_t^2)K_tX_t^2+2M_t\Theta_tX_t^2\right]dt\\
&+\left(2C_tK_t+\Theta_t\right)X_t^2dW_t,
\end{align*}
where $\tilde{\Phi}$ related to $K$ and $\Theta$, which is a rewriting of $\Phi$ (as in equation (\ref{4.3})), is given by
\begin{align*}
\tilde{\Phi}(t,x)=&-(R_t+D_t^2K_t)^{-1}(B_tK_tx+D_t\Theta_tx+C_tD_tK_tx)^2\notag\\
&\qquad+(2A_tK_t+2C_t\Theta_t+S_t+C_t^2K_t)x^2\\
=&\inf_u\Big[\big(R_t+D_t^2K_t\big) u^2\notag\\
&\qquad +2(B_tK_tx+D_t\Theta_tx+C_tD_tK_tx)u\notag\\
&\qquad+(2A_tK_t+2C_t\Theta_t+S_t+C_t^2K_t)x^2\Big].
\end{align*}
Then 
\begin{align*}
d\left(K_t^2X_t^4\right)=&d\left(V(t,X_t)\right)^2\\
=&2K_tX_t^2dV(t,X_t)+d\langle V(t,X_t),V(t,X_t)\rangle_2\\
=&-2K_t\tilde{\Phi}(t,X_t)X_t^2dt+\left[\left(2C_tK_t+\Theta_t\right)X_t^2\right]^2dt\\
&+\left[\left(2A_t+C_t^2\right)K_t+2C_t\Theta_t\right]K_tX_t^4dt\\
&+2\left(2C_tK_t+\Theta_t\right)K_tX_t^4dW_t.
\end{align*}

Taking the integration and using the fact $\frac{1}{2}a^2\le (a+b)^2+b^2$, we have
\begin{align*}
\frac{1}{2}\int_0^T\Theta_t^2X_t^4dt\le &\int_0^T\left[\left(2C_tK_t+\Theta_t\right)X_t^2\right]^2dt+ \int_0^T4C_t^2K_t^2X_t^4dt\\
=&G^2X_T^4-x^4+2\int_0^TK_t\tilde{\Phi}(t,X_t)X_t^2dt\\
&-\int_0^T\left[\left(2A_t+C_t^2\right)K_t+2C_t\Theta_t\right]K_tX_t^4dt\\
&+2\int_0^T(2C_tK_t+\Theta_t)K_tX_t^4dW_t\\
&+\int_0^T4C_t^2K_t^2X_t^4dt\\
\le&G^2X_T^4+2\int_0^T[(2A_t+C_t^2)K_t+2C_t\Theta_t+S_t]K_tX_t^4dt\\
&-\int_0^T\left[\left(2A_t+C_t^2\right)K_t+2C_t\Theta_t\right]K_tX_t^4dt\\
&+2\int_0^T(2C_tK_t+\Theta_t)K_tX_t^4dW_t\\
&+\int_0^T4C_t^2K_t^2X_t^4dt\,. 
\end{align*}
Since the left side is positive, we have
\begin{align*}
\mathbb{E}\left(\int_0^T\Theta_t^2X_t^4dt\right)^2\le &L\mathbb{E}[G^4X_T^8]+\frac{1}{4}\mathbb{E}\left(\int_0^T\Theta_t^2X_t^4dt\right)^2\\
&+L\mathbb{E}\left(\int_0^T\left[(2A_t+C_t^2)K_t^2+S_tK_t\right]X_t^4dt\right)^2\\
&+L\mathbb{E}\int_0^T(2C_tK_t+\Theta_t)^2K_t^2X_t^8dt\\
\le&L\mathbb{E}[G^4X_T^8]+\frac{1}{4}\mathbb{E}\left(\int_0^T\Theta_t^2X_t^4dt\right)^2\\
&+L\mathbb{E}\int_0^T\left[(2A_t+C_t^2)^2K_t^4+S_t^2K_t^2\right]X_t^8dt\\
&+L\mathbb{E}\int_0^TC_t^2K_t^4X_t^8dt\\
&+\frac{1}{4}\mathbb{E}\left(\int_0^T\Theta_t^2X_t^4dt\right)^2+L\mathbb{E}\left[\sup_{t\in[0,T]}K_t^2X_t^4\right]^2,
\end{align*}
where $L$, which may vary from line to line, denotes a positive constant only depend on $T$. Recall the process $\mathcal{Y}$ defined by  equation (\ref{mathcaly}), we have
\begin{align*}
&\mathbb{E}\left(\int_0^T\Theta_t^2\mathcal{Y}_{0,t}^4dt\right)^2\\
\le& L\mathbb{E}\left[\left(\sup_{t\in[0,T]}\mathcal{Y}_{0,t}^8\right)\times\left(G^4+\int_0^T\left[A_t^2K_t^4+C_t^2K_t^4+C_t^4K_t^4+S_t^2K_t^2\right]dt\right)\right]\\
\le&L\mathbb{E}\left[\left(\sup_{t\in[0,T]}\mathcal{Y}_{0,t}^{16}\right)+\left(\sup_{t\in[0,T]}K_t^{16}\right)+\left(\sup_{t\in[0,T]}K_t^{8}\right)+G^8\right]\\
&\quad + L\mathbb{E}\int_0^T\left(A_t^8+C_t^8+C_t^{16}+S_t^8\right)dt\\
<&+\infty.
\end{align*}
Finally, we conclude
\begin{align*}
\mathbb{E}\int_0^T\Theta_t^2dt\le&\mathbb{E}\left[\left(\sup_{t\in[0,T]}\mathcal{Z}_{0,t}^4\right)\times\int_0^T\Theta_t^2\mathcal{Y}_{0,t}^4dt\right]\\
\le&\left[\mathbb{E}\left(\sup_{t\in[0,T]}\mathcal{Z}_{0,t}^8\right)\right]^{\frac{1}{2}}\times\left[\mathbb{E}\left(\int_0^T\Theta_t^2\mathcal{Y}_{0,t}^4dt\right)^2\right]^{\frac{1}{2}}\\
<&+\infty.
\end{align*}
Here $\mathcal{Z}$ is given by
\begin{align}\label{mathcalz}
\mathcal{Z}_{s,t}=\mathcal{Y}_{s,t}^{-1}=exp\left\{-\int_s^t(A_r-\frac{1}{2}C_r^2)dr-\int_s^tC_rdW_r-\int_s^tF_r\circ dB_r^H\right\}.
\end{align}
\end{proof}

Up to now, combining with Lemma \ref{lem4.3}
and Theorem \ref{the5.1}, we prove the existence of the solution to the BSRE (\ref{bsre}). Then we give the uniqueness result.

\begin{theorem}
There is the unique solution pair $(K,\Theta)$ to BSRE (\ref{bsre}).
\end{theorem}
\begin{proof}
Let both $(K,\Theta)$ and $(\tilde{K},\tilde{\Theta})$ be the solutions to the BSRE (\ref{bsre}). Denote
\begin{align*}
\hat{\Phi}(t;K,\Theta)=&-(R_t+D_t^2K)^{-1}(B_tK+D_t\Theta+C_tD_tK)^2\\
&+(2A_t+C_t^2)K+2C_t\Theta+S_t\\
=&\inf_u\hat{G}(t,K,\Theta,u),
\end{align*}
with 
\begin{align*}
\hat{G}(t,K,\Theta,u)=&\big(R_t+D_t^2K\big) u^2\notag\\
& +2(B_tK+D_t\Theta+C_tD_tK)u\notag\\
&+(2A_tK+2C_t\Theta+S_t+C_t^2K).
\end{align*}
So that, using the notation $\Delta K=\tilde{K}-K, \,\Delta\Theta=\tilde{\Theta}-\Theta$, we have
\begin{align*}
-d\Delta K_t=&\left[\hat{\Phi}(t;\tilde{K}_t,\tilde{\Theta}_t)-\hat{\Phi}(t;K_t,\Theta_t)\right]dt+F_t\Delta K_t\circ dW_t^H-\Delta\Theta_tdW_t.
\end{align*}
Applying   It$\hat{\rm o}$'s formula to $(\Delta K_t)^2$, we have
\begin{align}\label{5.2}
\mathbb{E}\left[(\Delta K_T)^2-(\Delta K_0)^2\right]=&-2\mathbb{E}\int_0^T \Delta K_t \left[\hat{\Phi}(t;\tilde{K}_t,\tilde{\Theta}_t)-\hat{\Phi}(t;K_t,\Theta_t)\right]dt\notag\\
&-2\mathbb{E}\int_0^TF_t\left(\Delta K_t\right)^2\circ dW_t^H+\mathbb{E}\int_0^T(\Delta \Theta_t)^2dt\notag\\
\ge&-2\mathbb{E}\int_0^T|\Delta K_t|\max_{i\in\{1,2\}}|\hat{G}(t,\tilde{K}_t,\tilde{\Theta}_t,u^i_t)-\hat{G}(t,K_t,\Theta_t,u^i_t)|dt\notag\\
&-2\mathbb{E}\int_0^TF_t\left(\Delta K_t\right)^2\circ dW_t^H+\mathbb{E}\int_0^T(\Delta \Theta_t)^2dt.
\end{align}
According to Theorem \ref{the4.5}, 
\begin{align*}
K_t=\inf_{u\in\mathcal{U}[t,T]}J(u;t,1),
\end{align*}
which means the uniqueness of $K$, since the lower greatest bound of $J(u;t,1)$ is unique. Then the uniqueness of $\Theta$ follows from inequality (\ref{5.2}).
\end{proof}

\section{Conclusion and discussion}
This paper fills  a fundamental gap in the literature on stochastic control problems driven by fractional Brownian motions by developing a stochastic dynamic programming principle for such non-Markovian systems. Despite the long-standing success of the DPP in the classical Brownian (Markovian)  setting, its extension to fractional noise has remained largely open due to the loss of the Markov property. By exploiting the structure of fractional Brownian motion and its underlying Brownian representation, we establish a rigorous dynamic programming framework and derive a stochastic HJB equation in the form of a backward stochastic partial differential equation driven jointly by fractional Brownian motion and standard Brownian motion. A verification theorem is further obtained, providing a precise characterization of optimal controls and clarifying the relationship between the dynamic programming approach and the stochastic maximum principle.

As an application, we study an  LQ optimal control problem of controlled systems driven by both Brownian motions and fractional Brownian motions. The value function is shown to admit a quadratic representation, and the optimal control is characterized explicitly through a backward stochastic Riccati equation. In contrast to the classical Markovian case, the resulting optimal control is adapted but generally non-Markovian, reflecting the intrinsic memory effect of fractional noise. The existence and uniqueness of solutions to the associated backward stochastic Riccati equation are also established, thereby completing the dynamic programming analysis for the linear–quadratic framework.

Notice that throughout this paper the diffusion coefficient $\sigma$ is assumed to be independent of the control. On the other hand, from the viewpoint of the stochastic maximum principle, Han et al. derived necessary optimality conditions for controlled systems driven by fractional Brownian motion (see Theorem 5.7 of \cite{han2013maximum}). Yet their conditions include the constraint $\sigma^*_u(t)p_t=0$, which essentially excludes any nontrivial control dependence in $\sigma$. Therefore, both the DPP approach and the maximum principle face intrinsic difficulties in this setting. Consequently, the optimal control problem for fractional Brownian motion–driven systems with control-dependent diffusion coefficients is still an open problem.

Here we present  a simple example to show there may not exist an optimal control in the control dependent diffusion case, which suggests that we may need to impose significantly more conditions or more restricted admissible control space. Let the state equation be
\begin{align}
X_t=\int_0^tu_t\circ dW_t^H,\quad t\in[0,T],\label{e.6.1}
\end{align}
with the cost functional
\begin{align}
J(u)=\mathbb{E}\left[\int_0^Tu_t^2dt-X_T\right].\label{e.6.2}
\end{align}
Through formula (\ref{2.4}), we have
\begin{align*}
J(u)=\mathbb{E}\int_0^Tu_t^2-\mathbb{D}_t^H[u_t]dt.
\end{align*}
Let 
\begin{align*}
u_t=\mathbf{1}_{[\delta,T]}(t)\cdot\delta^k\cdot(W_t^H-W_{t-\delta}^H),\quad \delta\in(0,T).
\end{align*}
Then, we have
\begin{align*}
\mathbb{E}|u_t|^2=\mathbf{1}_{[\delta,T]}(t)\cdot\delta^{2k+2H},
\end{align*}
and 
\begin{align*}
\mathbb{D}_t^H[u_t]=&H(2H-1)\int_0^t|t-s|^{2H-2}D_s^H[u_t]ds\\
=&\mathbf{1}_{[\delta,T]}(t)H(2H-1)\delta^k\int_{t-\delta}^t|t-s|^{2H-2}ds\\
=&\mathbf{1}_{[\delta,T]}(t)H(2H-1)^2\delta^{k+2H-1}.
\end{align*}
So that
\begin{align*}
J(u)=(T-\delta)\left[\delta^{2k+2H}-H(2H-1)^2\delta^{k+2H-1}\right].
\end{align*}
Choose $-H<k<1-2H$ and let $\delta\to 0^+$, it follows
\begin{align*}
    J(u)\to-\infty,
\end{align*}
which means there doesn't exist an optimal control
for the problem \eqref{e.6.1}-\eqref{e.6.2}.

\bibliographystyle{siamplain}
\bibliography{ex_supplement}
\end{document}